\documentclass{amsart}
\usepackage{times,hyperref,tikz}
\usepackage[english]{babel}

\renewcommand{\normalsize}{\fontsize{10}{13pt}\selectfont}
\renewcommand{\small}{\fontsize{9}{11pt}\selectfont}

\renewcommand{\large}{\fontsize{12}{14pt}\selectfont}

\def\proc#1{\medbreak\noindent{\it #1}\hspace{1ex}\ignorespaces}
\def\ep{\noindent{\hfill $\Box$}}

\newtheorem{theo}{Theorem}[section]
\newtheorem{prop}[theo]{Proposition}
\newtheorem{lemma}[theo]{Lemma}
\newtheorem{cor}[theo]{Corollary}
\numberwithin{equation}{section}

\title[Dynamical properties of the negative beta transformation]{\large Dynamical properties of the negative beta transformation}

\author[Lingmin Liao]{\normalsize Lingmin Liao}
\address{LAMA  CNRS UMR 8050,
Universit\'e Paris-Est - Cr\'eteil - Val-de-Marne,  UFR Sciences et
Technologie, 61 Avenue du G\'en\'eral de Gaulle, 94010 Cr\'eteil
Cedex, France} 
\email{lingmin.liao@u-pec.fr}

\author{Wolfgang Steiner}
\address{LIAFA, CNRS, Universit\'e Paris Diderot -- Paris 7, Case 7014, 75205 Paris Cedex 13, France}
\email{steiner@liafa.jussieu.fr}

\begin{document}
\maketitle

\begin{abstract}
We analyse dynamical properties of the negative beta transformation, which has been studied recently by Ito and Sadahiro.
Contrary to the classical beta transformation, the density of the absolutely continuous invariant measure of the negative beta transformation may be zero on certain intervals. 
By investigating this property in detail, we~prove that the $(-\beta)$-transformation is exact for all $\beta>1$, confirming a conjecture of G\'ora, and intrinsic, which completes a study of Faller. 
We also show that the limit behaviour of the $(-\beta)$-expansion of $1$ when $\beta$ tends to~$1$ is related to the Thue-Morse sequence. 

A consequence of the exactness is that every Yrrap number, which is a $\beta>1$ such that the $(-\beta)$-expansion of $1$ is eventually periodic, is a Perron number.
This extends a well-known property of Parry numbers.
However, the set of Parry numbers is different from the set of Yrrap numbers.
\end{abstract}

\section{Introduction}
For a real number $\beta>1$, the \emph{$\beta$-transformation} is defined by
\[
T_\beta:\, [0,1) \to [0,1)\,,\quad x \mapsto \beta x - \lfloor \beta x \rfloor\,.
\]
R\'enyi~\cite{Renyi57} was the first to use it for representing real numbers in base~$\beta$, generalising expansions in integer bases.
The (greedy) \emph{$\beta$-expansion} of $x \in [0,1)$~is
\[
x = \frac{\lfloor\beta x\rfloor}{\beta} + \frac{\lfloor\beta\,T_\beta(x)\rfloor}{\beta^2} + \frac{\lfloor\beta\,T_\beta^2(x)\rfloor}{\beta^3} + \cdots\,.
\]

In this paper, we study the \emph{$(-\beta)$-transformation} (negative beta transformation)
\[
T_{-\beta}:\,(0,1] \to (0,1]\,,\quad x \mapsto -\beta x + \lfloor
\beta x\rfloor + 1\,.
\]
Note that $T_{-\beta}(x) = -\beta x - \lfloor -\beta x \rfloor$ except for finitely many points, hence $T_{-\beta}$ is a natural modification of the $\beta$-transformation, avoiding the discontinuity of $x \mapsto -\beta x - \lfloor -\beta x \rfloor$ at~$0$.
For $x \in (0,1]$, set $d_{-\beta}(x) = d_{-\beta,1}(x)\, d_{-\beta,2}(x) \cdots$ with
\[
d_{-\beta,1}(x) = \lfloor \beta x\rfloor + 1\,, \quad
d_{-\beta,n}(x) = d_{-\beta,1}(T_{-\beta}^{n-1}(x)) \quad \mbox{for}\
n \ge 1\,,
\]
then the \emph{$(-\beta)$-expansion} of~$x$ is
\[
x  = \sum_{k=1}^\infty \frac{-d_{-\beta,k}(x)}{(-\beta)^k} = \frac{d_{-\beta,1}(x)}{\beta} - \frac{d_{-\beta,2}(x)}{\beta^2} + \frac{d_{-\beta,3}(x)}{\beta^3} - \frac{d_{-\beta,4}(x)}{\beta^4} + \cdots\,.
\]
Examples of a $\beta$-transformation and a $(-\beta)$-transformation are depicted in Figure~\ref{f:1}.

\begin{figure}[ht]
\centering
\begin{tikzpicture}[scale=4]
\draw(0,0)--(1,0)--(1,1)--(0,1)--cycle (.618,0)--(.618,1);
\draw[thick](0,0)--(.618,1) (.618,0)--(1,.618);

\begin{scope}[shift={(1.5,0)}]
\draw(0,0)--(1,0)--(1,1)--(0,1)--cycle (.618,0)--(.618,1);
\draw[thick](0,1)--(.618,0) (.618,1)--(1,.382);
\end{scope}
\end{tikzpicture}
\caption{\mbox{$\beta$-transformation (left) and $(-\beta)$-transformation (right),~$\beta = \frac{1+\sqrt{5}}{2}$.}}
\label{f:1}
\end{figure}

The map $T_{-\beta}$ has been studied recently by G\'ora~\cite{Gora07} and Faller~\cite{Faller08}.
Ito and Sadahiro~\cite{Ito-Sadahiro09} considered a $(-\beta)$-transformation on
the interval $\big[\frac{-\beta}{\beta+1}, \frac{1}{\beta+1}\big)$, defined by
\[
\textstyle x \mapsto -\beta x - \big\lfloor \frac{\beta}{\beta+1} - \beta x\big\rfloor\,.
\]
We remark that their transformation is conjugate to our $T_{-\beta}$ through the conjugacy function $\phi(x)= \frac{1}{\beta+1}-x$.
So, all the results in~\cite{Ito-Sadahiro09} can be translated into our
case directly, in particular their $(-\beta)$-expansion of $\phi(x) \in \big[\frac{-\beta}{\beta+1}, \frac{1}{\beta+1}\big)$ is
\[
\phi(x) = \sum_{k=1}^\infty \frac{d_{-\beta,k}(x)-1}{(-\beta)^k}\,.
\]

Both $\beta$-transformation and $(-\beta)$-transformation are
examples of piecewise monotonic transformations. 
Dynamical properties of general piecewise monotonic transformations have
been investigated in the literature
\cite{Hofbauer79,Hofbauer81a,Hofbauer81b,Keller78,Wagner79}.
However, for the specific case of the $(-\beta)$-transformation, the
detailed dynamical properties are unknown. 
As pointed out in~\cite{Ito-Sadahiro09}, for any $\beta>1$, $T_{-\beta}$~admits a
unique absolutely continuous invariant measure (hence ergodic), with the density
\[
h_{-\beta}(x) = \sum_{n\geq 1, \ T_{-\beta}^{n}(1) \geq x }
\frac{1}{(-\beta)^n}\,.
\]
Contrary to the classical $ \beta$-transformation, Ito and Sadahiro~\cite{Ito-Sadahiro09} gave an example such that the density $h_{-\beta}$ is zero on some intervals.

In the present paper, we fully study this phenomenon. An interval
on which the density equals zero will be called a \emph{gap}. From a
result of Keller~\cite{Keller78}, one can deduce that, for fixed
$\beta>1$, the number of gaps is finite.
In Theorem~\ref{t:gaps}, we describe the set of the
gaps and show that, as $\beta$ decreases, the number of gaps
forms an increasing sequence
\[
0 < 1 < 2 < 5 < 10 < \cdots < \lfloor 2^{n+1}/3 \rfloor < \cdots\,.
\]
Figures~\ref{f:2} and~\ref{f:3} show examples with $2$ gaps and $5$ gaps respectively.
The endpoints of the gaps are determined by the orbit of~$1$,
which is described by the sequence~$d_{-\beta}(1)$. 

We also show (Theorem~\ref{t:limit}) that
\[
\lim_{\beta\to1} d_{-\beta}(1) = 211222112112112221122 \cdots\,,
\]
which is the fixed point of the following morphism:
\[
\varphi:\, 1 \mapsto 2\,, \quad 2 \mapsto 211\,.
\]
Here, $\varphi$ is a morphism on finite words on the alphabet $\{1,2\}$, which is naturally extended to the set of infinite words $\{1,2\}^{\mathbb{N}}$ as in~\cite{Lothaire83}, and $\lim_{\beta\to1} d_{-\beta}(1) = u$ means that longer and longer prefixes of $d_{-\beta}(1)$ agree with prefixes of $u$ when~$\beta \to 1$.
We remark that this sequence has been studied e.g.\ in \cite{AABBJPS95,Dubickas06,Dubickas07}, and is known to be the smallest aperiodic sequence in $\{1,2\}^{\mathbb{N}}$ with the property that all its proper suffixes are smaller than itself with respect to the alternate order.
Moreover, by adding a $1$ at the beginning, we obtain the sequence of run lengths in the Thue-Morse sequence.
It is shown in~\cite{Ito-Sadahiro09} that every sequence $d_{-\beta}(1)$, $\beta > 1$, has the property that all its proper suffixes are smaller than itself with respect to the alternate order.

A dynamical system $(X, T)$ is said to be \emph{locally eventually onto} or
topologically exact if, for any non-empty open subset $U\subset X$,
there exists a positive integer $n$ such that $f^n(U)=X$.
The main theorem of this paper (Theorem \ref{t:main-mixing}) asserts
that $T_{-\beta}$ is locally eventually onto on $(0,1] \setminus
G(\beta)$, where $G(\beta)$ is the (finite) union of the gaps.

The notion of topological exactness is derived from the exactness in
ergodic theory introduced by Rohlin~\cite{Rohlin61}. 
A~measure-preserving transformation $T$ on a probability space $(X, \mathcal{B}, \mu)$ is \emph{exact} if and only if, for any positive measure
subset $A$ with $T^n(A)\in \mathcal{B}$ $(n\geq 0)$, one has
$\lim_{n\to\infty} \mu(T^n(A)) =1$ (see~\cite[p.~125]{Pollicott-Yuri98}). As a
corollary of our main theorem, we confirm G\'ora's conjecture~\cite{Gora07} that all
$(-\beta)$-transformations are exact with respect to the unique absolutely continuous invariant measure.

An invariant measure of a dynamical system is called a \emph{maximal
entropy measure} if it maximises entropy. If there is a unique
maximal entropy measure, the dynamical system is called
\emph{intrinsic}. As another application of our main theorem, we give a
complete answer about the uniqueness of maximal entropy
measures discussed in the thesis of Faller~\cite{Faller08}. We prove that all
$(-\beta)$-transformations are intrinsic (Corollary~\ref{c:Faller}).

Finally, we use the main result to show that every number $\beta > 1$ with eventually periodic $T_{-\beta}$-orbit of~$1$ is a Perron number.
Recall that a number $\beta > 1$ with the corresponding property for $T_\beta$ is called a \emph{Parry number} in reference to~\cite{Parry60}, and note that $T_{-\beta}(x) = 1 - T_\beta(x)$ for all $x \in (0,1)$.
Because of this orientation reversing property, we call a number $\beta > 1$ with eventually periodic $T_{-\beta}$-orbit of~$1$ an \emph{Yrrap number}.
Parry numbers are also known to be Perron numbers. 
However, not every Perron number is a Parry number or an Yrrap number.
In Section~\ref{sec:proofs-main-results}, we give examples showing that the set of Parry numbers and the set of Yrrap numbers do not include each other.

The paper is organised as follows. 
In Section~\ref{Sec:Results}, we state the main results. 
In Section~\ref{sec:polyn-expans}, we establish properties of the sequences~$d_{-\beta}(1)$, which allow describing the structure of the gaps in Section~\ref{sec:structure-gaps}.
In Section~\ref{sec:locally-event-onto}, we show that $T_{-\beta}$ is locally eventually onto on $(0,1] \setminus G(\beta)$.
The proofs of the main results are completed in Section~\ref{sec:proofs-main-results}.

We remark that, independently of our work, some of our results have been proved for slightly more general transformations by Hofbauer~\cite{Hofbauer}.

\section{Main results}\label{Sec:Results}
For each $n \ge 1$, let $\gamma_n$ be the (unique) positive real number defined by
\[
\gamma_n^{g_n+1} = \gamma_n + 1\,, \quad \mbox{with} \quad g_n = \lfloor 2^{n+1}/3 \rfloor\,,
\]
and set $\gamma_0 = \infty$.
Then
\[
2 > \gamma_1 > \gamma_2 > \gamma_3 > \cdots > 1\,.
\]
Note that $\gamma_1$ is the golden ratio and that $\gamma_2$ is the
smallest Pisot number.

For each $n \ge 0$ and $1 < \beta < \gamma_n$, set
\[
\textstyle \mathcal{G}_n(\beta) = \big\{G_{m,k}(\beta) \mid 0 \le m < n,\, 0 \le k < \frac{2^{m+1}+(-1)^m}{3}\big\}\,,
\]
with open intervals
\[
G_{m,k}(\beta) = \left\{\begin{array}{cl}\big(T_{-\beta}^{2^{m+1}+k}(1),\, T_{-\beta}^{(2^{m+2}-(-1)^m)/3+k}(1)\big) & \mbox{if $k$ is even}, \\[1ex] \big(T_{-\beta}^{(2^{m+2}-(-1)^m)/3+k}(1),\, T_{-\beta}^{2^{m+1}+k}(1)\big) & \mbox{if $k$ is odd},\end{array}\right.
\]
in particular $\mathcal{G}_0(\beta)$ is the empty set.

\begin{theo}\label{t:gaps}
For any $\gamma_{n+1} \leq \beta < \gamma_n$, $n \ge 0$, the set of gaps of the transformation~$T_{-\beta}$ is~$\mathcal{G}_n(\beta)$, which consists of $g_n = \lfloor 2^{n+1}/3 \rfloor$ disjoint non-empty intervals.
\end{theo}

We define therefore
\[
G(\beta) = \bigcup_{I\in\mathcal{G}_n(\beta)} I \quad \mbox{if}\ \gamma_{n+1} \leq \beta < \gamma_n,\, n \ge 0\,.
\]

\begin{theo}\label{t:main-mixing}
For any $\beta > 1$, $T_{-\beta}$ is locally eventually onto on $(0,1] \setminus G(\beta)$,
\[
T_{-\beta}^{-1}(G(\beta)) \subset G(\beta) \quad \text{and} \quad \lim_{k\to\infty} \lambda\big(T_{-\beta}^{-k}(G(\beta))\big) = 0\,.
\]
\end{theo}

G\'ora~\cite{Gora07} proved that, for $\beta > \gamma_2$, the
transformation $T_{-\beta}$ is exact, and he conjectured that this
would hold for all $\beta > 1$. We confirm his conjecture.

\begin{cor}\label{c:Gora}
For any $\beta > 1$, the transformation $T_{-\beta}$ is exact with
respect to its unique absolutely continuous invariant measure.
\end{cor}

\proc{Proof.}
The unicity of the absolutely continuous invariant measure follows from~\cite{Li-Yorke78}, as already noted in~\cite{Ito-Sadahiro09}.
Now, the exactness is a direct consequence of Theorem~\ref{t:main-mixing}.
\ep\medbreak

Faller~\cite{Faller08} proved that, for $\beta> \sqrt[3]{2}$, $T_{-\beta}$ admits a unique maximal entropy measure
and left the case $\beta\leq \sqrt[3]{2}$. We give a complete
answer.

\begin{cor}\label{c:Faller}
For any $\beta > 1$, the transformation $T_{-\beta}$ has a unique
maximal entropy measure, hence is intrinsic.
\end{cor}

\proc{Proof.}
By Theorem~\ref{t:main-mixing}, the transformation $T_{-\beta}$ is locally eventually onto on $(0,1] \setminus G(\beta)$.
A~result of Walters \cite[Theorem~16, p.~140]{Walters78} completes the proof.
\ep\medbreak

The key to understanding the $(-\beta)$-transformation is to know the $(-\beta)$-expansion of~$1$.
Here, the main tool will be the morphism $\varphi:\, 1 \mapsto 2,\ 2 \mapsto 211$.
The expansion of~$1$ under $T_{-\gamma_n}$ can be described by this morphism, and the expansion of~$1$ under $T_{-\beta}$ tends to the fixed point of the morphism when $\beta$ tends to~$1$ (from above).

\begin{theo}\label{t:limit}
For every $n \ge 1$, we have
\[
d_{-\gamma_n}(1) = \varphi^{n-1}\big(2\,\overline{1}\big)\,,
\]
where $\overline{1} = 111\cdots$.
If $1 < \beta \le \gamma_n$, then $d_{-\beta}(1)$ starts with $\varphi^n(2)$, hence
\[
\lim_{\beta\to1} d_{-\beta}(1) = \lim_{n\to\infty} \varphi^n\big(2\, \overline{1}\big) = 211222112112112221122\cdots\,.
\]
\end{theo}

Mas\'akov\'a and Pelantov\'a~\cite{Masakova-Pelantova} showed that the Yrrap numbers (called Ito--Sadahiro numbers by them) are algebraic integers with all conjugates having modulus less than~$2$, thus all Yrrap numbers $\beta \ge 2$ are Perron numbers.
Again, we give a complete answer.

\begin{theo} \label{t:perron}
Every Yrrap number is a Perron number.
\end{theo}

At the end of this section, we recall some identities.
Note that
\begin{align}
|\varphi^n(1)| & = \frac{2^{n+1}+(-1)^n}{3} = g_n + \frac{1+(-1)^n}{2}\,, \label{id-1} \\
|\varphi^n(2)| & = |\varphi^{n+1}(1)| = \frac{2^{n+2}-(-1)^n}{3} = g_{n+1} + \frac{1-(-1)^n}{2}\,, \label{id-2}\\
|\varphi^n(11)| & = 2\cdot|\varphi^{n}(1)| = \frac{2^{n+2}+2\cdot(-1)^n}{3} = g_{n+1} + \frac{1+(-1)^n}{2}\,, \label{id-3}
\end{align}
where $|w|$ denotes the length of the word~$|w|$.
In particular, $|\varphi^n(1)|$ and $|\varphi^n(2)|$ are odd, $|\varphi^n(21)| = 2^{n+1}$, $g_n$ is odd if and only if $n$ is odd, and we can write
\begin{align*}
\mathcal{G}_n(\beta) & = \big\{G_{m,k}(\beta) \mid 0 \le m < n,\, 0 \le k < |\varphi^m(1)|\big\}\,, \\[1ex]
G_{m,k}(\beta) & = \left\{\begin{array}{cl}\big(T_{-\beta}^{|\varphi^m(21)|+k}(1),\, T_{-\beta}^{|\varphi^m(2)|+k}(1)\big) & \mbox{if $k$ is even}, \\[1ex] \big(T_{-\beta}^{|\varphi^m(2)|+k}(1),\, T_{-\beta}^{|\varphi^m(21)|+k}(1)\big) & \mbox{if $k$ is odd}.\end{array}\right.
\end{align*}

\section{Polynomials and expansions} \label{sec:polyn-expans}

To study the trajectories of $1$ under~$T_{-\beta}$, we define the maps
\[
f_{\beta,a}:\, \mathbb{R} \to \mathbb{R}\,,\ x \mapsto -\beta x + a\,, \quad a \in \{1,2\}\,.
\]
The composition of maps $f_{\beta,a_1}, \ldots, f_{\beta,a_k}$, is denoted by
\[
f_{\beta,a_1 \cdots a_k} = f_{\beta,a_k} \circ \cdots \circ f_{\beta,a_1}\,.
\]
Since
\[
f_{\beta,a_1 \cdots a_k}(1) = (-\beta)^k + \sum_{j=1}^k a_k (-\beta)^{k-j}\,,
\]
we define the polynomial
\[
P_{a_1 \cdots a_k} = (-X)^k + \sum_{j=1}^k a_k (-X)^{k-j} \in \mathbb{Z}[X]
\]
for every word $a_1 \cdots a_k \in \{1,2\}^k$, $k \ge 1$. 
Then $f_{\beta,a_1 \cdots a_k}(1) = P_{a_1 \cdots a_k}(\beta)$, thus
\[
T_{-\beta}^k(1) = f_{\beta,d_{-\beta,1}(1) \cdots d_{-\beta,k}(1)}(1) = P_{d_{-\beta,1}(1) \cdots d_{-\beta,k}(1)}(\beta)\,.
\]

For the proof of Theorem~\ref{t:limit}, we use the following polynomial identities.

\begin{lemma} \label{l:poly}
For $1 \le j < k$, we have
\[
P_{a_1 \cdots a_k} = (-X)^{k-j} \big(P_{a_1 \cdots a_j} - 1\big) + P_{a_{j+1} \cdots a_k}\,.
\]
\end{lemma}

\proc{Proof.}
We can deduce the identity directly from the definition.
\ep\medbreak

\begin{lemma} \label{l:poly2}
For every $n \ge 0$, we have
\[
X^{\frac{1+(-1)^n}{2}} P_{\varphi^n(2)} + X^{\frac{1-(-1)^n}{2}} P_{\varphi^n(11)} = X + 1 = X^{\frac{1+(-1)^n}{2}} + X^{\frac{1-(-1)^n}{2}}\,.
\]
\end{lemma}

\proc{Proof.}
The second equation holds for even and odd~$n$, thus for all $n \in \mathbb{Z}$.
Since 
\[
X P_2 + P_{11}= X (-X + 2) + X^2 - X + 1\,,
\]
the first equation holds for $n=0$.
For every $n \ge 0$, we have
\begin{align*}
X^{\frac{1+(-1)^{n+1}}{2}} P_{\varphi^{n+1}(2)} + X^{\frac{1-(-1)^{n+1}}{2}} P_{\varphi^{n+1}(11)} & = X^{\frac{1-(-1)^n}{2}} P_{\varphi^n(211)} + X^{\frac{1+(-1)^n}{2}} P_{\varphi^n(22)}.
\end{align*}
Then, using that $|\varphi^n(11)|$ is even, $|\varphi^n(2)|$ is odd and Lemma~\ref{l:poly}, we obtain
\begin{align*}
X^{\frac{1-(-1)^n}{2}} P_{\varphi^n(211)} + X^{\frac{1+(-1)^n}{2}} P_{\varphi^n(22)}
& = X^{\frac{1-(-1)^n}{2}+|\varphi^n(11)|} \big(P_{\varphi^n(2)} - 1\big) \\
& \hspace{-11em} + X^{\frac{1-(-1)^n}{2}} P_{\varphi^n(11)} - X^{\frac{1+(-1)^n}{2}+|\varphi^n(2)|} \big(P_{\varphi^n(2)} - 1\big) + X^{\frac{1+(-1)^n}{2}} P_{\varphi^n(2)} \\
& = X^{\frac{1+(-1)^n}{2}} P_{\varphi^n(2)} + X^{\frac{1-(-1)^n}{2}} P_{\varphi^n(11)}\,.
\end{align*}
Therefore, the lemma follows inductively.
\ep\medbreak

The following lemma is a consequence of Lemma~\ref{l:poly2}.

\begin{lemma} \label{l:poly2a}
For every $n \ge 0$, the words $\varphi^n(2)$ and $\varphi^n(11)$ agree on the first $g_{n+1}-1$ letters and differ on the $g_{n+1}$-st letter.
\end{lemma}

\proc{Proof.}
The polynomials $X^{\frac{1+(-1)^n}{2}} P_{\varphi^n(2)}$ and $X^{\frac{1-(-1)^n}{2}} P_{\varphi^n(11)}$ have degree $g_{n+1} + 1$ by~\eqref{id-2} and~\eqref{id-3}.
Thus the coefficient of $X^{g_{n+1}+1-j}$ in their sum is equal to the difference of the $j$-th letter in $\varphi^n(2)$ and $\varphi^n(11)$, for $1 \le j \le \min(|\varphi^n(2)|, |\varphi^n(11)|) = g_{n+1}$.
Since the sum of these polynomials is equal to $X + 1$, by Lemma~\ref{l:poly2}, the first $g_{n+1}-1$ letters in $\varphi^n(2)$ and $\varphi^n(11)$ are equal, and the $g_{n+1}$-st letters differ.
\ep\medbreak

\begin{lemma} \label{l:poly3}
For every $n \ge 0$, we have
\begin{align}
1 - P_{\varphi^n(1)} & = X^{\frac{1+(-1)^n}{2}} \prod_{m=0}^{n-1} \big(X^{|\varphi^m(1)|} - 1\big)\,, \label{e:poly1} \\
P_{\varphi^n(21)}  - 1 & = \Big(X^{g_{n+1}+1} - X^{g_n+1} - X^{\frac{1+(-1)^n}{2}}\Big) \prod_{m=0}^{n-1} \big(X^{|\varphi^m(1)|} - 1\big)\,, \label{e:poly2} \\
P_{\varphi^n(21)} - P_{\varphi^n(2)} & = \big(X^{g_{n+1}+1} - X - 1\big) \prod_{m=0}^{n-1} \big(X^{|\varphi^m(1)|} - 1\big)\,, \label{e:poly3} \\
P_{\varphi^n(2)} - P_{\varphi^{n+1}(2)} & = \big(X^{g_{n+1}+1} - X - 1\big) \prod_{m=0}^n \big(X^{|\varphi^m(1)|} - 1\big)\,. \label{e:poly4}
\end{align}
\end{lemma}

\proc{Proof.}
Since $1 - P_1 = X$, \eqref{e:poly1} holds for $n=0$.
By Lemmas~\ref{l:poly2} and~\ref{l:poly}, we obtain that
\begin{align*}
1 - P_{\varphi^{n+1}(1)} & = 1 - P_{\varphi^n(2)} \\
& = X^{-(-1)^n} \big(P_{\varphi^n(11)} - 1\big) \\
& = X^{-(-1)^n} \big(-X^{|\varphi^n(1)|} \big(P_{\varphi^n(1)} - 1\big) + P_{\varphi^n(1)} - 1\big) \\
& = X^{-(-1)^n} \big(X^{|\varphi^n(1)|} - 1\big) \big(1 - P_{\varphi^n(1)}\big)
\end{align*}
for every $n \ge 0$, hence \eqref{e:poly1} holds inductively for all $n \ge 0$.

Applying first Lemma~\ref{l:poly}, then equation~\eqref{e:poly1} for $n$ and~$n+1$, we obtain that
\begin{align*}
P_{\varphi^n(21)} - 1 & = - X^{|\varphi^n(1)|} \big(P_{\varphi^n(2)} - 1\big) + P_{\varphi^n(1)} - 1 \\
& = \Big(X^{|\varphi^n(1)|+\frac{1-(-1)^n}{2}} \big(X^{|\varphi^n(1)|} - 1\big) - X^{\frac{1+(-1)^n}{2}}\Big) \prod_{m=0}^{n-1} \big(X^{|\varphi^m(1)|} - 1\big)\,.
\end{align*}
By \eqref{id-1} and~\eqref{id-3}, we have $|\varphi^n(1)|+\frac{1-(-1)^n}{2} = g_n + 1$ and $2 \cdot |\varphi^n(1)|+\frac{1-(-1)^n}{2} = g_{n+1} + 1$, which gives~\eqref{e:poly2}.

Combining \eqref{e:poly1} for~$n+1$ and \eqref{e:poly2}, we obtain~\eqref{e:poly3}.

Finally, \eqref{e:poly4} follows from \eqref{e:poly1} for $n+1$ and~$n+2$.
\ep\medbreak

With the help of Lemma~\ref{l:poly3}, we will show the following proposition, where $\eta_n > 1$, $n \ge 1$, is defined by
\[
\eta_n^{g_n+1} = \eta_n^{g_{n-1}+1} + \eta_n^{\frac{1-(-1)^n}{2}}\,.
\]
It can be easily verified that
\[
2 = \eta_1 > \gamma_2^2 > \gamma_1 > \eta_2 > \gamma_2 > \eta_3 > \cdots > 1\,.
\]

\begin{prop} \label{p:expansion}
Let $\beta > 1$, $n \ge 1$. 
Then $d_{-\beta}(1)$ starts with $\varphi^n(2)$ if and only if $\beta < \gamma_2^2$ in case $n = 1$, $\beta \le \eta_n$ in case $n$ is even, $\beta < \eta_n$ in case $n \ge 3$ is odd.
\end{prop}

For the proof of the proposition, we use the following lemmas.
For the sake of readability, we often omit the dependence on $\beta$ in $T_{-\beta}$, $d_{-\beta}$ and $f_{\beta,a}$ in the sequel.

\begin{lemma}\label{l:first-digit}
Let $x\in (0,1]$ and $1 < \beta < 2$.
Then
\begin{itemize}
\itemsep0pt
\item
$d_1(x) = 2$ if and only if $f_2(x) \le 1$,
\item
$d_1(x) = 1$ if and only if $f_1(x) > 0$,
\item
$d_1(x)\, d_2(x) = 11$ if and only if $0 < f_{11}(x) < 1$.
\end{itemize}
\end{lemma}

\proc{Proof.}
For $a \in \{1,2\}$, we have $d_1(x) = a$ if and only if $f_a(x) \in (0,1]$.
The inequalities $f_2(x) > - \beta + 2 > 0$ and $f_1(x) < 1$ prove the first two points.
Noting that $f_1(x) < 1$ is equivalent to $x > 0$, the third point follows from the second point.
\ep\medbreak

\begin{lemma}\label{l:digit-joint}
Let $\beta > 1$. 
If $d(1)$ starts with $a_1 \cdots a_k$ and $d(f_{a_1 \cdots a_k}(1))$ starts with $b_1\cdots b_j$, then $d(1)$ starts with $a_1 \cdots a_k\, b_1\cdots b_j$.
\end{lemma}

\proc{Proof.}
If $d_1(1) \cdots d_k(1) = a_1 \cdots a_k$, then $T^k(1) = f_{a_1 \cdots a_k}(1)$, and the lemma follows from $d_{k+1}(1) \cdots d_{k+j}(1) = d_1(T^k(1)) \cdots d_j(T^k(1))$.
\ep\medbreak

\begin{lemma} \label{l:phin1}
Let $x \in (0,1]$, $n \ge 1$, such that $d(x)$ starts with the first $g_n - 1$ letters of~$\varphi^n(1)$.
Then $d(x)$ starts with~$\varphi^n(1)$ if and only if $f_{\varphi^n(1)}(x) \le 1$ in case $n$ is odd, $0 < f_{\varphi^n(1)}(x) < 1$ in case $n$ is even.
\end{lemma}

\proc{Proof.}
Assume that $d(x)$ starts with the first $g_n - 1$ letters of~$\varphi^n(1)$, $n \ge 1$.
If $n$ is odd, then $g_n = |\varphi^n(1)|$, and $\varphi^n(1)$ ends with~$2$, thus $f_2\, T^{|\varphi^n(1)|-1}(1) = f_{\varphi^n(1)}(1)$.
By Lemma~\ref{l:first-digit}, we have $d_{|\varphi^n(1)|}(1) = d_1(T^{|\varphi^n(1)|-1}(1)) = 2$ if and only if $f_2\, T^{|\varphi^n(1)|-1}(1) \le 1$.
If $n$ is even, then $g_n = |\varphi^n(1)|-1$, and $\varphi^n(1)$ ends with~$11$, thus $f_{11}\, T^{|\varphi^n(1)|-2}(1) = f_{\varphi^n(1)} (1)$.
By Lemma~\ref{l:first-digit}, $d_{|\varphi^n(1)|-1}(1)\, d_{|\varphi^n(1)|}(1) = 11$ is equivalent to $0 < f_{11}\, T^{|\varphi^n(1)|-2}(1) < 1$.
\ep\medbreak

\begin{lemma} \label{l:continuous}
Let $\beta > 1$, $n \ge 1$.
If $d(1)$ starts with~$\varphi^n(2)$, then $T^{g_n-1}$ is continuous on $\big[T^{|\varphi^n(1)|}(1),\, 1\big]$, and $1/\beta$ is in the interior of $T^{g_n-1}\big(\big[T^{|\varphi^n(1)|}(1),\, 1\big]\big)$.
\end{lemma}

\proc{Proof.}
Let $\beta > 1$ and $n \ge 1$ such that $d(1)$ starts with $\varphi^n(2) = \varphi^{n-1}(211)$, then $d(T^{|\varphi^n(1)|}(1))$ starts with~$\varphi^{n-1}(11)$. 
By Lemma~\ref{l:poly2a}, the words $\varphi^{n-1}(2)$ and $\varphi^{n-1}(11)$ share the first $g_n-1$ letters and differ on the $g_n$-th letter.
This proves that $T^{g_n-1}$ is continuous on $\big[T^{|\varphi^n(1)|}(1),\, 1\big]$, and that $1/\beta$ is in the interior or the right endpoint of $T^{g_n-1}\big(\big[T^{|\varphi^n(1)|}(1),\, 1\big]\big)$.
If $T^{g_n-1}(1) = 1/\beta$ with odd $n$, then $T^{|\varphi^n(1)|}(1) = 1$. 
If $T^{|\varphi^n(1)|+g_n-1}(1) = T^{g_{n+1}-1}(1) = 1/\beta$ with even~$n$, then $T^{|\varphi^n(2)|}(1) = 1$.
Both situations are impossible since $f_{\varphi^n(1)}(1) < 1$ for all $\beta > 1$, $n \ge 0$, by~\eqref{e:poly1}.
Therefore, $1/\beta$ is in the interior of $T^{g_n-1}\big(\big[T^{|\varphi^n(1)|}(1),\, 1\big]\big)$.
\ep\medbreak

\begin{lemma} \label{l:ordering}
For any $1 < \beta \le \eta_n$, $n \ge 1$, we have $f_{\varphi^{n-1}(1)}(1) \le f_{\varphi^n(2)}(1) < 1$. 
\end{lemma}

\proc{Proof.}
By~\eqref{e:poly1}, we have $f_{\varphi^n(1)}(1) < 1$ for all $\beta > 1$, $n \ge 0$.

For all $1 < \beta \le \eta_n$, $n \ge 1$, we have $f_{\varphi^{n-1}(21)}(1) \le 1$ by~\eqref{e:poly2}.
Since $f_{\varphi^{n-1}(1)}$ is order reversing, we get that $f_{\varphi^n(2)}(1) = f_{\varphi^{n-1}(211)}(1) \ge f_{\varphi^{n-1}(1)}(1)$.
\ep\medbreak

\proc{Proof of Proposition~\ref{p:expansion}.}
Let $\beta > 1$. 
First note that $d_1(1) = 2$ if and only if~$\beta < 2$. 
Then, by Lemma~\ref{l:first-digit}, $d(1)$ starts with $\varphi(2)$ if and only if $0 < f_{\varphi(2)}(1) < 1$.
By~\eqref{e:poly1}, we have $1 - f_{\varphi(2)}(1) = \beta (\beta - 1)^2 > 0$, and $1 - f_{\varphi(2)}(1) < 1$ if and only if $\beta - 1 < \beta^{-1/2}$, i.e., $\beta < \gamma_2^2$.
This proves the case $n = 1$ of the proposition.

Suppose now that $d(1)$ starts with $\varphi^{n+1}(2) = \varphi^n(211)$ for some $n \ge 1$.
Then we must have $f_{\varphi^n(21)}(1) \le 1$, which is equivalent to $\beta \le \eta_{n+1}$ by~\eqref{e:poly2}.
Suppose moreover that the proposition is already shown for~$n$.
Since $\eta_{n+1} < \min(\eta_n, \gamma_2^2)$, this implies that $d(1)$ starts with~$\varphi^n(2)$.
Furthermore, $\eta_{n+1} < \gamma_n$ and \eqref{e:poly4} yield that $f_{\varphi^n(2)}(1) > f_{\varphi^n(1)}(1)$. 
By Lemma~\ref{l:continuous}, the first $g_n-1$ letters of $d(x)$ are equal to those of $d(1)$ for all $x \in \big[f_{\varphi^n(1)}(1),\, 1\big]$.
Therefore, Lemma~\ref{l:phin1} implies that $d(1)$ starts with $\varphi^n(21)$ if and only if $f_{\varphi^n(21)}(1) \le 1$ in case $n$ is odd, $0 < f_{\varphi^n(21)}(1) < 1$ in case $n$ is even.
For even~$n$, we must therefore have $\beta < \eta_{n+1}$. 
Note that $f_{\varphi^n(21)}(1) > 0$ is guaranteed by $f_{\varphi^n(21)}(1) > f_{\varphi^n(1)}(1) = T^{|\varphi^n(1)|}(1)$, which follows from $f_{\varphi^n(2)}(1) < 1$ and the order reversing property of~$f_{\varphi^n(1)}$.
By Lemma~\ref{l:ordering}, we have $f_{\varphi^n(1)}(1) \le f_{\varphi^n(211)}(1) < 1$, thus $f_{\varphi^n(21)}(1) \in \big[f_{\varphi^n(1)}(1),\, 1\big]$, which implies that $d(f_{\varphi^n(21)}(1))$ starts with~$\varphi^n(1)$.
Putting everything together, we obtain that $d(1)$ starts with $\varphi^n(211) = \varphi^{n+1}(2)$ if and only if $\beta < \eta_{n+1}$ with even~$n$, or $\beta \le \eta_{n+1}$ with odd~$n$.
\ep\medbreak

We conclude the section with the proof of Theorem~\ref{t:limit}.

\proc{Proof of Theorem~\ref{t:limit}.}
By Proposition~\ref{p:expansion}, $d_{-\beta}(1)$ starts with $\varphi^n(2) $ for any $1 < \beta \le \gamma_n$, $n \ge 1$. 
By~\eqref{e:poly3}, we have $f_{\gamma_n,\varphi^{n-1}(21)}(1) = f_{\gamma_n,\varphi^{n-1}(2)}(1)$, thus $T_{-\gamma_n}^{|\varphi^{n-1}(21)|}(1) = T_{-\gamma_n}^{|\varphi^{n-1}(2)|}(1)$, i.e., $T_{-\gamma_n}^{|\varphi^{n-1}(2)|}(1) $ is a fixed point of $T_{-\gamma_n}^{|\varphi^{n-1}(1)|}$. 
Hence, by Lemma \ref{l:digit-joint}, $d_{-\gamma_n}(1) = \varphi^{n-1}(2\,\overline{1})$. 
Since $\lim_{n\to\infty} \gamma_n = 1$, we obtain that $d_{-\beta}(1)$ starts with $\varphi^n(2)$ for larger and larger~$n$ when $\beta \to 1$, thus $\lim_{\beta\to1} d_{-\beta}(1) = \lim_{n\to\infty} \varphi^n(2\, \overline{1})$.
\ep\medbreak

\section{Structure of the gaps} \label{sec:structure-gaps}

\begin{figure}[ht]
\centering
\begin{tikzpicture}[scale=16.5]
\small
\fill[black!10](.0531,0)--(.0724,0)--(.0724,.0531)--(.7,.0531)--(.7,.0724)--(.0724,.0724)--(.0724,.7)--(.0531,.7)--(.0531,.0724)--(0,.0724)--(0,.0531)--(.0531,.0531);
\fill[black!10](.1563,0)--(.25,0)--(.25,.1563)--(.7,.1563)--(.7,.25)--(.25,.25)--(.25,.7)--(.1563,.7)--(.1563,.25)--(0,.25)--(0,.1563)--(.1563,.1563);
\fill[black!10](.3245,0)--(.3417,0)--(.3417,.3245)--(.7,.3245)--(.7,.3417)--(.3417,.3417)--(.3417,.7)--(.3245,.7)--(.3245,.3417)--(0,.3417)--(0,.3245)--(.3245,.3245);
\fill[black!10](.4478,0)--(.5242,0)--(.5242,.4478)--(.7,.4478)--(.7,.5242)--(.5242,.5242)--(.5242,.7)--(.4478,.7)--(.4478,.5242)--(0,.5242)--(0,.4478)--(.4478,.4478);
\fill[black!10](.6185,0)--(.6338,0)--(.6338,.6185)--(.7,.6185)--(.7,.6338)--(.6338,.6338)--(.6338,.7)--(.6185,.7)--(.6185,.6338)--(0,.6338)--(0,.6185)--(.6185,.6185);
\draw[dotted]
(0,.25)node[left]{$t_1$}--(.7,.25) (.25,0)node[below]{$t_1$}--(.25,.7)
(0,.1563)node[left]{$t_2$}--(.7,.1563) (.1563,0)node[below]{$t_2$}--(.1563,.7)
(0,.5242)node[left]{$t_3$}--(.7,.5242) (.5242,0)node[below]{$t_3$}--(.5242,.7)
(0,.4478)node[left]{$t_4$}--(.7,.4478) (.4478,0)node[below]{$t_4$}--(.4478,.7)
(0,.6338)node[left]{$t_5$}--(.7,.6338) (.6338,0)node[below]{$t_5$}--(.6338,.7)
(0,.3245)node[left]{$t_6$}--(.7,.3245) (.3245,0)node[below]{$t_6$}--(.3245,.7)
(0,.0724)node[left]{$t_7$}--(.7,.0724) (.0724,0)node[below]{$t_7$}--(.0724,.7)
(0,.6185)node[left]{$t_8$}--(.7,.6185) (.6185,0)node[below]{$t_8$}--(.6185,.7)
(0,.3417)node[left]{$t_9$}--(.7,.3417) (.3417,0)node[below]{$t_9$}--(.3417,.7)
(0,.0531)node[left]{$t_{10}$}--(.7,.0531) (.0531,0)node[below]{$t_{10}$}--(.0531,.7)
(.1563,0)--(.25,0) (.4478,0)--(.5242,0) (.7,.1563)--(.7,.25) (.7,.4478)--(.7,.5242) (.5242,.7)--(.4478,.7) (.25,.7)--(.1563,.7) (0,.5242)--(0,.4478) (0,.25)--(0,.1563) (.3889,.1563)--(.3889,.25) (.3889,.4478)--(.3889,.5242);
\draw(0,.1563)--(0,0)node[left]{$0$}node[below]{$0$}--(.1563,0) (.25,0)--(.4478,0) (.5242,0)--(.7,0)node[below]{$t_0$}--(.7,.1563) (.7,.25)--(.7,.4478) (.7,.5242)--(.7,.7)--(.5242,.7) (.4478,.7)--(.25,.7) (.1563,.7)--(0,.7)node[left]{$t_0$}--(0,.5242) (0,.4478)--(0,.25) (.3889,0)--(.3889,.1563) (.3889,.25)--(.3889,.4478) (.3889,.5242)--(.3889,.7);
\draw[thick](0,.7)--(.1563,.5242) (.25,.1563)--(.389,0) (.389,.7)--(.4478,.6338) (.5242,.4478)--(.7,.25);
\draw[thick,dotted](.1563,.5242)--(.25,.1563) (.4478,.6338)--(.5242,.4478);
\begin{scope}[shift={(0,-.035)}]
\draw(0,0)--node[below]{$F_{3,5}$}(.0531,0)--node[below]{$G_{3,2}$}(.0724,0)--node[below]{$F_{3,2}$}(.1563,0) (.25,0)--node[below]{$F_{3,1}$}(.3245,0)--node[below]{$G_{3,1}$}(.3417,0)--node[below]{$F_{3,4}$}(.4478,0) (.5242,0)--node[below]{$F_{3,3}$}(.6185,0)--node[below]{$G_{3,0}$}(.6338,0)--node[below]{$F_{3,0}$}(.7,0)
(0,-.04)--(0,.005) (.0531,-.005)--(.0531,.005) (.0724,-.005)--(.0724,.005) (.1563,-.04)--(.1563,.005) (.25,-.04)--(.25,.005) (.3245,-.005)--(.3245,.005) (.3417,-.005)--(.3417,.005) (.4478,-.04)--(.4478,.005) (.5242,-.04)--(.5242,.005) (.6185,-.005)--(.6185,.005) (.6338,-.005)--(.6338,.005) (.7,-.04)--(.7,.005)
(0,-.035)--node[below]{$F_{2,2}$}(.1563,-.035) (.25,-.035)--node[below]{$F_{2,1}$}(.4478,-.035) (.5242,-.035)--node[below]{$F_{2,0}$}(.7,-.035);
\draw[dotted](.1563,0)--node[below]{$G_{1,0}$}(.25,0) (.4478,0)--node[below]{$G_{2,0}$}(.5242,0);
\end{scope}
\end{tikzpicture}
\caption{The $(-\beta)$-transformation for $\beta = 9/8$ with gaps $G_{m,k}$, support of the invariant measure $F_{3,0} \cup \cdots \cup F_{3,5}$, and $t_k = T_{-\beta}^k(1)$.
Ratios between $G_{1,0}$, $G_{2,0}$ and the rest of the picture are not respected.
In reality, $G_{1,0}$ is almost $10$ times the size of~$G_{2,0}$, and $G_{2,0}$ is almost $5$ times the size of~$F_{2,0}$.}
\label{f:3}
\end{figure}
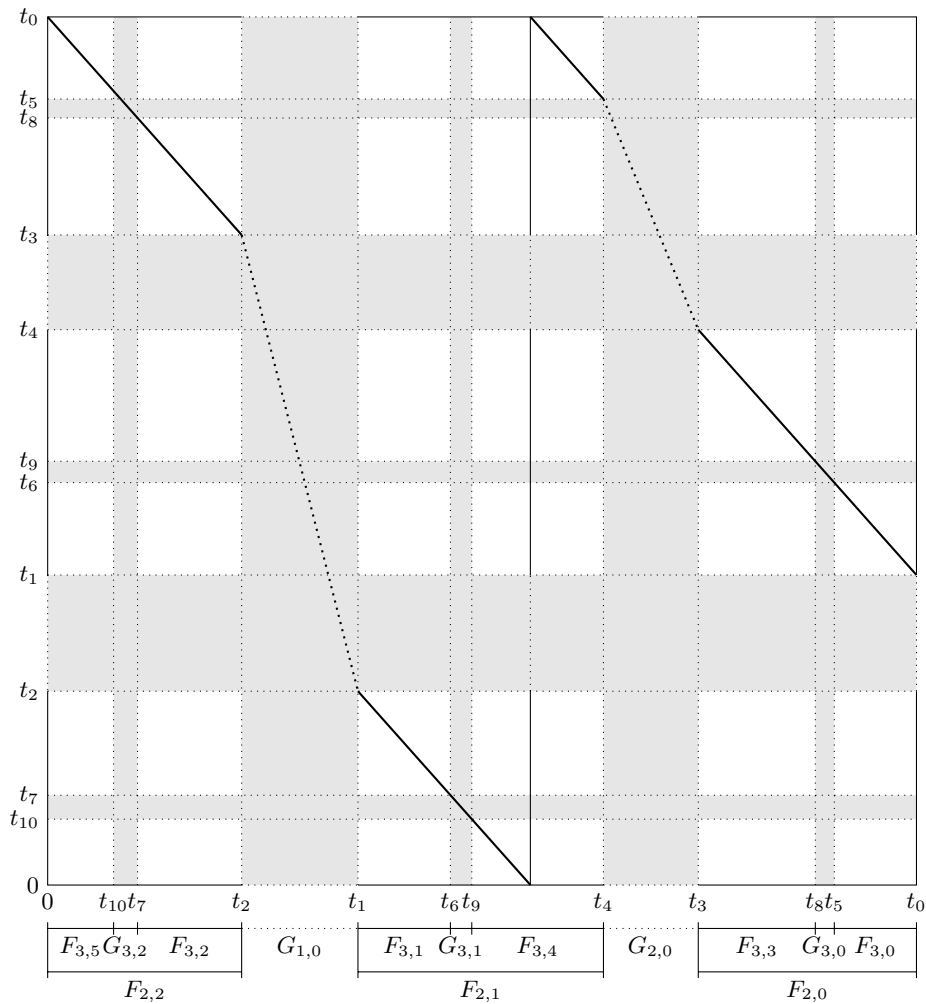

We will show that the support of the invariant measure is
\[
F(\beta) = (0,1] \setminus G(\beta)\,.
\]
For $\beta > 1$, $n \ge 0$, let
\[
\mathcal{F}_n(\beta) = \big\{F_{n,k}(\beta) \mid 0 \le k \le g_n\big\}
\]
with
\begin{align*}
F_{n,k}(\beta) & = \left\{\begin{array}{cl}\big[T_{-\beta}^{|\varphi^n(1)|+k}(1),\, T_{-\beta}^k(1)\big] & \mbox{if $k$ is even},\ k < g_n, \\[1ex] \big[T_{-\beta}^k(1),\, T_{-\beta}^{|\varphi^n(1)|+k}(1)\big] & \mbox{if $k$ is odd},\ k < g_n,\end{array}\right. \\
F_{n,g_n}(\beta) & =  \left\{\begin{array}{cl}\big(0,\, T_{-\beta}^{|\varphi^n(1)|-1}(1)\big] & \mbox{if $n$ is even,} \\[1ex] \big(0,\, T_{-\beta}^{|\varphi^n(2)|-1}(1)\big] & \mbox{if $n$ is odd.}\end{array}\right.
\end{align*}
(We will consider these sets only in the case when the left number is smaller than the right number.)
Figures~\ref{f:2} and~\ref{f:3} show examples of the decomposition of $(0,1]$ into these sets.

\begin{figure}[ht]
\centering
\begin{tikzpicture}[scale=10]
\small
\fill[black!10](.0625,0)--(.75,0)--(.75,.0625)--(1,.0625)--(1,.75)--(.75,.75)--(.75,1)--(.0625,1)--(.0625,.75)--(0,.75)--(0,.0625)--(.0625,.0625);
\fill[black!10](.84765625,0)--(.921875,0)--(.921875,.84765625)--(1,.84765625)--(1,.921875)--(.921875,.921875)--(.921875,1)--(.84765625,1)--(.84765625,.921875)--(0,.921875)--(0,.84765625)--(.84765625,.84765625);
\draw[dotted](0,.75)node[left]{$t_1$}--(1,.75) (.75,0)node[below]{$t_1$}--(.75,1)
(0,.0625)node[left]{$t_2$}--(1,.0625) (.0625,0)node[below]{$t_2$}--(.0625,1)
(0,.921875)node[left]{$t_3$}--(1,.921875) (.921875,0)node[below]{$t_3$}--(.921875,1)
(0,.84765625)node[left]{$t_4$}--(1,.84765625) (.84765625,0)node[below]{$t_4$}--(.84765625,1);
\draw(0,0)node[left]{$0$}node[below]{$0$}--(1,0)node[below]{$t_0$}--(1,1)--(0,1)node[left]{$t_0$}--cycle (.8,0)--(.8,1);
\draw[thick](0,1)--(.8,0) (.8,1)--(1,.75);
\begin{scope}[shift={(0,-.07)}]
\draw(0,0)--node[below]{$F_{2,2}$}(.0625,0)--node[below]{$G_{1,0}$}(.75,0)--node[below]{$F_{2,1}$}(.84765625,0)--node[below]{$G_{2,0}$}(.921875,0) --node[below]{$F_{2,0}$}(1,0);
\draw(0,-.08)--(0,.01) (.0625,-.08)--(.0625,.01) (.75,-.08)--(.75,.01) (.84765625,-.01)--(.84765625,.01) (.921875,-.01)--(.921875,.01) (1,-.08)--(1,.01);
\draw(0,-.07)--node[below]{$F_{1,1}$}(.0625,-.07) (.75,-.07)--node[below]{$F_{1,0}$}(1,-.07);
\end{scope}
\end{tikzpicture}
\caption{The $(-\beta)$-transformation for $\beta = 5/4$, with set of gaps $G_{1,0} \cup G_{2,0}$ and support of the invariant measure $F_{2,2} \cup F_{2,1} \cup F_{2,0}$.
Here, $t_k = T_{-\beta}^k(1)$.}
\label{f:2}
\end{figure}
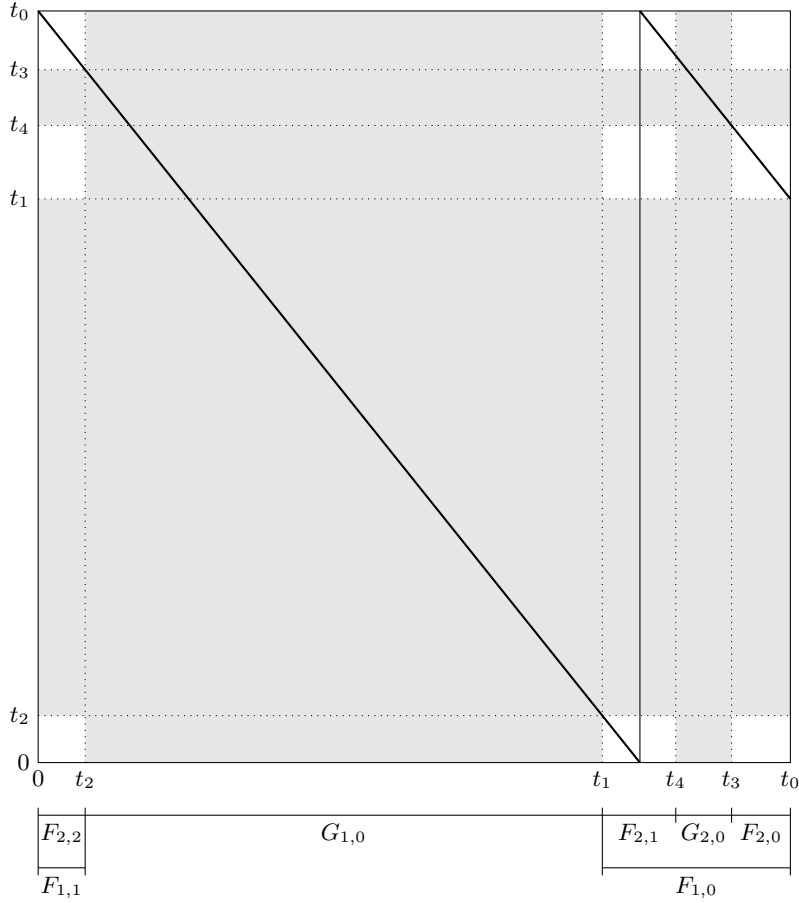

As in Section~\ref{sec:polyn-expans}, we often omit the dependence on $\beta$ in the following.

\begin{lemma}
Let $n \ge 2$ and $1 < \beta \le \eta_n$ if $n$ is even, $1 < \beta < \eta_n$ if $n$ is odd; or let $n = 1$, $1 < \beta < \gamma_2^2$.
Then we have
\renewcommand{\theenumi}{\roman{enumi}}
\begin{enumerate}
\item \label{enum3}
$F_{n,k} = T^k(F_{n,0})$ for all $0 \le k < g_n$,
\item \label{enum4}
$1/\beta$ is in the interior of $F_{n,g_n-1}$,
\item \label{enum5}
$T^{g_n}(F_{n,0}) = F_{n,g_n} \cup F_{n,0}$, $T(F_{n,g_n}) = F_{n+1,0}\setminus \{1\}$, if $n$ is odd,
\item \label{enum6}
$T^{g_n}(F_{n,0}) = F_{n,g_n} \cup  F_{n+1,0}$, $T(F_{n,g_n}) = F_{n,0}\setminus \{1\}$, if $n$ is even.
\end{enumerate}
\end{lemma}

\proc{Proof.}
Let $\beta$ and $n$ have the properties of the statement of the lemma.
Then Proposition~\ref{p:expansion} and Lemma~\ref{l:continuous} give the points~\ref{enum3}.\ and~\ref{enum4}.

If $n$ is odd, then $F_{n,g_n-1} = \big[T^{g_{n+1}-1}(1),\, T^{g_n-1}(1)\big]$, thus 
\begin{align*}
T(F_{n,g_n-1}) & = \big(0,\, T^{g_{n+1}}(1)\big] \cup \big[T^{g_n}(1),\, 1\big]  = \big(0,\, T^{|\varphi^n(2)|-1}(1)\big] \cup \big[T^{|\varphi^n(1)|}(1),\, 1\big] \\
& = F_{n,g_n} \cup F_{n,0}\,.
\end{align*}
Since $d(1)$ starts with $\varphi^n(2)$ and $\varphi^n(2)$ ends with~$1$, we have
\[
T(F_{n,g_n}) = T\big(\big(0,\, T^{|\varphi^n(2)|-1}(1)\big]\big) = \big[T^{|\varphi^n(2)|}(1), 1\big) = F_{n+1,0} \setminus \{1\}\,,
\]
i.e., \ref{enum5}.\ holds.
If $n$ is even, then $F_{n,g_n-1} = \big[T^{g_n-1}(1),\, T^{g_{n+1}-1}(1)\big]$ and 
\begin{align*}
T(F_{n,g_n-1}) & = \big(0, T^{g_n}(1)\big] \cup \big[T^{g_{n+1}}(1),\, 1\big] = \big(0, T^{|\varphi^n(1)|-1}(1)\big] \cup \big[T^{|\varphi^n(2)|}(1), 1\big] \\
& = F_{n,g_n} \cup  F_{n+1,0}\,. 
\end{align*}
Now, $\varphi^n(1)$ ends with~$1$, which gives $T(F_{n,g_n}) = F_{n,0} \setminus \{1\}$, i.e., \ref{enum6}.\ holds.
\ep\medbreak

\begin{prop} \label{p:partition}
Let $1 < \beta < \gamma_n$, $n \ge 0$. 
Then the elements of $\mathcal{F}_n$ and $\mathcal{G}_n$ are intervals of positive length which form a partition of $(0,1]$.
If $n \ge 1$, then
\renewcommand{\theenumi}{\roman{enumi}}
\begin{enumerate}
\itemsep3pt
\item \label{enum1}
$F_{n-1,k} = F_{n,|\varphi^{n-1}(1)|+k} \cup G_{n-1,k} \cup F_{n,k}$ for all $0 \le k < |\varphi^{n-1}(1)|$, 
\item \label{enum7}
$G_{m,k} = T^k(G_{m,0})$ for all $0 \le m < n$, $0 \le k < |\varphi^m(1)|$,
\item \label{enum8}
$T^{|\varphi^m(1)|}(G_{m,0}) = G_{m,0} \cup F_{m,0} \setminus F_{m+1,0} \supset G_{m,0}$ for all $0 \le m < n$.
\end{enumerate}
\end{prop}

\proc{Proof.}
We have $\mathcal{F}_0 = \{F_{0,0}\} = \{(0,1]\}$ and $\mathcal{G}_0 = \emptyset$, thus the elements of $\mathcal{F}_0$ and $\mathcal{G}_0$ form the trivial partition of $(0,1]$.
Let now $1 < \beta < \gamma_n$, $n \ge 1$, and assume that $\mathcal{F}_{n-1} \cup \mathcal{G}_{n-1}$ forms a partition of $(0,1]$.
Lemma~\ref{l:poly3} and the order reversing property of $f_{\varphi^{n-1}(1)}$ give
\[
f_{\varphi^{n-1}(1)}(1) < f_{\varphi^{n-1}(21)}(1) < f_{\varphi^{n-1}(2)}(1) < f_{\varphi^n(2)}(1) < 1\,.
\]
For $n = 1$, we have $f_{21}(1) > 0$. 
Therefore, $F_{n-1,0}$ splits into $F_{n,|\varphi^{n-1}(1)|}$, $G_{n-1,0}$ and $F_{n,0}$, which are intervals of positive length.
This shows~\ref{enum1}.\ in case $n = 1$.
By Proposition~\ref{p:expansion} and Lemma~\ref{l:continuous}, $T^{g_{n-1}-1}$ is continuous on $F_{n-1,0}$, thus every $F_{n-1,k}$, $0 \le k < g_{n-1}$, also splits into $F_{n,|\varphi^{n-1}(1)|+k}$, $G_{n-1,k}$ and $F_{n,k}$. 
If $n$ is even, then this proves~\ref{enum1}.\ since $g_{n-1} = |\varphi^{n-1}(1)|$. 
In this case, we have $F_{n-1,g_{n-1}} = F_{n,g_n}$, thus $\mathcal{F}_n \cup \mathcal{G}_n$ is a refinement of the partition $\mathcal{F}_{n-1} \cup \mathcal{G}_{n-1}$.
If $n$ is odd, then we use that $T^{|\varphi^{n-1}(1)|}$ is continuous on $G_{n-1,0} \cup F_{n-1,0}$ since both $d(T^{|\varphi^{n-1}(21)|}(1))$ and $d(1)$ start with~$\varphi^{n-1}(1)$.
In this case, $F_{n-1,g_{n-1}} = F_{n,|\varphi^{n-1}(1)|-1}$ splits into $F_{n,g_n}$, $G_{n-1,|\varphi^{n-1}(1)|-1}$ and $F_{n,|\varphi^{n-1}(1)|-1}$, thus~\ref{enum1}.\ holds for odd $n$ too, and $\mathcal{F}_n \cup \mathcal{G}_n$ is again a refinement of the partition $\mathcal{F}_{n-1} \cup \mathcal{G}_{n-1}$. 

For all $0 \le m < n$, the continuity of $T^{|\varphi^m(1)|}$ on $G_{m,0}$ gives~\ref{enum7}.
Moreover, $T^{|\varphi^m(1)|}(G_{m,0}) = \big(T^{|\varphi^m(21)|}(1),\, T^{|\varphi^{m+1}(2)|}(1)\big)$ and $F_{m+1,0} \subset F_{m,0}$ imply~\ref{enum8}. 
\ep\medbreak

\begin{cor} \label{c:TFG}
Let $1 < \beta < \gamma_n$, $n \ge 0$.
Then
\[
\textstyle T_{-\beta}\big(\bigcup_{I\in\mathcal{F}_n(\beta)} I\big) = \bigcup_{I\in\mathcal{F}_n(\beta)} I \quad \mbox{and} \quad T_{-\beta}^{-1}\big(\bigcup_{I\in\mathcal{G}_n(\beta)} I\big) \subseteq \bigcup_{I\in\mathcal{G}_n(\beta)} I\,.
\]
If $\gamma_{n+1} \le \beta < \gamma_n$, then $F(\beta)  = \bigcup_{I\in\mathcal{F}_n(\beta)} I$ (and $G(\beta)  = \bigcup_{I\in\mathcal{G}_n(\beta)} I$ by definition).
\end{cor}

\begin{prop} \label{p:gap}
Let $\beta > 1$ and $\mu$ be an invariant measure of $T_{-\beta}$ which is absolutely continuous with respect to the Lebesgue measure.
Then $\mu(G(\beta)) = 0$.
\end{prop}

\proc{Proof.}
If $\beta \ge \gamma_1$, then $\mu(G(\beta)) = \mu(\emptyset) = 0$.
We will show, by induction on~$n$, that
\begin{equation} \label{e:muGnk}
\mu(G_{n,k}) = 0 \quad \mbox{for all}\ n \ge 0\,,\ 0 \le k < |\varphi^n(1)|\,,\ 1 < \beta < \gamma_{n+1}\,.
\end{equation}
Since $G = \bigcup_{m=0}^{n-1} \bigcup_{k=0}^{|\varphi^m(1)|-1} G_{m,k}$ for $\gamma_{n+2} \le \beta < \gamma_{n+1}$, this proves the proposition.

For $1 < \beta < \gamma_1$, we have $G_{0,0} = \big(T^2(1),\, T(1)\big) = \big(f_1\, T(1),\, T(1)\big)$ and 
\[
T^{-1}(G_{0,0}) = f_1^{-1}(G_{0,0}) = \big(f_1^{-1}\, T(1),\, T(1)\big) \subset G_{0,0}\,.
\]
Since $\mu(T^{-1}(G_{0,0})) = \mu(G_{0,0})$, we obtain $\mu\big(\big(f_1\, T(1),\, f_1^{-1}\, T(1)\big]\big) = 0$.
Iteratively, we get, for all $k \ge 0$,
\[
\mu\big(\big(f_1^{-2k+1}\, T(1),\, f_1^{-2k-1}\, T(1)\big]\big) = 0\,, \quad \mu\big(\big[f_1^{-2k-2}\, T(1),\, f_1^{-2k}\, T(1)\big)\big) = 0\,.
\]
This gives $\mu\big(\big(T^2(1),\, \frac{1}{\beta+1}\big)\big) = 0$ and
$\mu\big(\big(\frac{1}{\beta+1},\, T(1)\big)\big) = 0$, since $f_1^{-1}$ is contracting with fixed point $\frac{1}{\beta+1}$. 
Then, the absolute continuity of $\mu$ implies $\mu(G_{0,0}) = 0$, hence \eqref{e:muGnk} holds for $n = 0$.

Now consider $1 < \beta < \gamma_{n+1}$, $n \ge 1$, and assume that \eqref{e:muGnk} holds for $n-1$.
By Proposition~\ref{p:partition}~\ref{enum1}.--\ref{enum8}., we have
\begin{equation} \label{e:TGmk}
T(G_{m,k-1}) = \left\{\begin{array}{ll}G_{m,k} & \mbox{if}\ 1 \le k < |\varphi^m(1)|,\, 0 \le m \le n, \\[1ex] G_{m,0} \cup F_{m+1,|\varphi^m(1)|} \cup G_{m+1,0} & \mbox{if}\ k = |\varphi^m(1)|,\, 0 \le m < n.\end{array}\right.
\end{equation}
Since $T^{-1}\big(\bigcup_{I\in\mathcal{G}_n} I\big) \subset \bigcup_{I\in\mathcal{G}_n} I$ by Corollary~\ref{c:TFG}, we obtain that
\[
T^{-1}(G_{n,0}) \subset G_{n,|\varphi^n(1)|-1} \cup G_{n-1,|\varphi^{n-1}(1)|-1} \ \mbox{and} \ T^{1-|\varphi^n(1)|}(G_{n,|\varphi^n(1)|-1}) = G_{n,0}\,.
\]
Thus, up to a set of $\mu$-measure zero, $T^{-|\varphi^n(1)|}(G_{n,0}) = f_{\varphi^n(1)}^{-1}(G_{n,0}) \subset G_{n,0}$, with
\begin{align*}
f_{\varphi^n(1)}^{-1}(G_{n,0}) & = f_{\varphi^n(1)}^{-1}
\big(\big(f_{\varphi^n(1)}\, T^{|\varphi^n(2)|}(1),\, T^{|\varphi^n(2)|}(1)\big)\big) \\
& = \big(f_{\varphi^n(1)}^{-1} T^{|\varphi^n(2)|}(1),\,
T^{|\varphi^n(2)|}(1)\big)\,.
\end{align*}
As in the case $n=0$, using the $T$-invariance of~$\mu$, we obtain that
\begin{align*}
\mu\big(\big(f_{\varphi^n(1)}^{-2k+1}\, T^{|\varphi^n(2)|}(1),\, f_{\varphi^n(1)}^{-2k-1}\, T^{|\varphi^n(2)|}(1)\big]\big) & = 0\,, \\
\mu\big(\big[f_{\varphi^n(1)}^{-2k-2}\, T^{|\varphi^n(2)|}(1),\, f_{\varphi^n(1)}^{-2k}\, T^{|\varphi^n(2)|}(1)\big)\big) & = 0\,,
\end{align*}
for all $k \ge 0$.
Thus $\mu\big(\big(T^{|\varphi^n(21)|}(1),\, y\big)\big) = 0$ and $\mu\big(\big(y,\, T^{|\varphi^n(2)|}(1)\big)\big) = 0$, where $y$ denotes the fixed point of $f_{\varphi^n(1)}$, i.e., $\mu(G_{n,0}) = 0$.
Since $\mu(G_{n,k}) = \mu(T^{-k}(G_{n,k})) = \mu(G_{n,0})$ for $1 \le k <  |\varphi^n(1)|$, \eqref{e:muGnk} holds for~$n$.
\ep\medbreak

\begin{prop} \label{p:gap2}
Let $\beta > 1$, then $\lim_{k\to\infty} \lambda\big(T_{-\beta}^{-k}(G(\beta))\big) = 0$.
\end{prop}

\proc{Proof.}
If $\beta \ge \gamma_1$, then $G(\beta) = \emptyset$ and the statement holds.
For $1 < \beta < \gamma_{n+1}$, $n \ge 0$, the preimage $T^{-k}(G_{n,0})$ can be written as a disjoint union
\[
T^{-k}(G_{n,0}) = \bigcup_{a_1\cdots a_k\in L_k} f_{a_1\cdots a_k}^{-1}(G_{n,0})
\]
with a set of words $L_k \subseteq \{1,2\}^k$.
Then we have $\mu(T^{-k}(G_{n,0})) = \#(L_k)\, \beta^{-k} \mu(G_{n,0})$.
By \eqref{e:TGmk} and Corollary~\ref{c:TFG}, $L_k$~consists of the length $k$ suffixes of all concatenations $1^{i_0} (\varphi(1))^{i_1} \cdots (\varphi^n(1))^{i_n}$, $i_0, i_1, \ldots, i_n \ge 0$.
Therefore, the number of elements in~$L_k$ is bounded by $k\choose{n}$, which grows polynomially in~$k$.
This yields $\lim_{k\to\infty} \mu(T^{-k}(G_{n,0})) = 0$.
As $T^{-j}(G_{n,j}) = G_{n,0}$ for $1 \le j <  |\varphi^n(1)|$, we also have $\lim_{k\to\infty} \mu(T^{-k}(G_{n,j})) = 0$, thus $\lim_{k\to\infty} \mu(T^{-k}(I)) = 0$ for all $I \in \mathcal{G}_{n+1}$.
\ep\medbreak

\section{Locally eventually onto} \label{sec:locally-event-onto}

\begin{lemma}\label{l:cover-F00}
Let $1 < \beta < \gamma_1$, $I$ be an interval of positive length in~$F(\beta)$, and $n \ge 1$ such that $\eta_{n+1} < \beta < \eta_n$ with odd $n \ge 1$, or $\eta_{n+1} \le \beta \le \eta_n$ with even $n \ge 2$.
Then there exists an $m \ge 0$ such that $T_{-\beta}^m(I) \supseteq F_{n,0}(\beta) \setminus \{1\}$.
\end{lemma}

\proc{Proof.}
We use ideas of G\'ora's proof of Proposition~8 in \cite{Gora07}.
We distinguish four cases.

\medskip\noindent
\textbf{Case~1:} odd~$n \ge 1$, $\eta_{n+1} < \beta < \gamma_n$.
In this case, we have $\beta^{2g_n} = \beta^{g_{n+1}} > \beta^{g_n} + \beta^{-1} \ge \beta + \beta^{-1} > 2$.
Therefore, the largest connected components of images of $I$ grow until they cover two consecutive discontinuities of $T^{2g_n}$ in the interior of~$F(\beta)$. 
Since the only discontinuity of~$T$ is $1/\beta \in F_{n,g_n-1}$, the images eventually also cover two consecutive discontinuities of $T^{2g_n}$ in the interior of $F_{n,g_n-1} = \big[T^{|\varphi^n(2)|-2}(1), T^{|\varphi^n(1)|-1}(1)\big]$.
We have
\begin{align*}
T^{2g_n}\big(\big[1/\beta, f_{2\varphi^n(1)}^{-1}(1)\big]\big) & = T^{g_n}\big(\big[1/\beta, T^{|\varphi^n(1)|-1}(1)\big]\big)\,, \\[1ex]
T^{2g_n}\big(\big[f_{11\varphi^n(1)}^{-1}(1), 1/\beta\big)\big) & = T^{g_n-1}\big(\big[1/\beta, T^{|\varphi^n(1)|-1}(1)\big)\big)\,, \\[1ex]
T^{2g_n}\big(\big(f_{2\varphi^n(1)}^{-1}(1), T^{|\varphi^n(1)|-1}(1)\big]\big) & = T^{g_n}\big(\big[T^{|\varphi^n(2)|-2}(1), 1/\beta\big)\big)\,, \\[1ex]
T^{2g_n}\big(\big[T^{|\varphi^n(2)|-2}(1), f_{11\varphi^n(1)}^{-1}(1)\big)\big) & = T^{g_n-1}\big(\big[T^{|\varphi^n(21)|-1}(1), 1/\beta\big)\big)\,.
\end{align*}
Using Lemma~\ref{l:continuous}, we obtain that $T^{2g_n}$ is continuous on these four intervals.
Note that $T^{|\varphi^n(2)|-2}(1) < f_{11\varphi^n(1)}^{-1}(1)$ since $f_{11\varphi^n(1)}\, T^{|\varphi^n(2)|-2}(1) = f_{\varphi^n(21)}(1) > 1$.
The discontinuities of $T^{2g_n}$ in the interior of $F_{n,g_n-1}$ are therefore $1/\beta$, $f_{2\varphi^n(1)}^{-1}(1)$ and $f_{11\varphi^n(1)}^{-1}(1)$.
With 
\[
T^{g_n+1}\big(\big[1/\beta, f_{2\varphi^n(1)}^{-1}(1)\big]\big) =F_{n,0} \quad \mbox{and} \quad T^{g_n+1}\big(\big[1/\beta, T^{|\varphi^n(1)|-1}(1)\big)\big) = F_{n,0} \setminus \{1\}\,,
\]
we obtain that $T^m(I) \supseteq F_{n,0} \setminus \{1\}$ for some $m \ge 0$.

\medskip\noindent
\textbf{Case~2:} odd~$n \ge 1$, $\gamma_{n+1} \le \beta \le \eta_{n+1}$. 
Since $\beta^{g_{n+1}+1} \ge \beta + 1 > 2$, we look at $T^{2g_n+1}$ on~$F_{n,g_n-1}$.
Note that $T^{2g_n}\big(\big[1/\beta, f_{2\varphi^n(1)}^{-1}(1)\big]\big) =F_{n,g_n-1}$, thus there is a discontinuity of $T^{2g_n+1}$ in the interior of $\big[1/\beta, f_{2\varphi^n(1)}^{-1}(1)\big]$.
Contrary to Case~1, $T^{2g_n+1}$ has no discontinuity in $F_{n,g_n-1}$ on the left of~$1/\beta$.
More precisely, 
\begin{align*}
T^{2g_n+1}\big(\big[1/\beta, f_{2\varphi^n(2)}^{-1}(1)\big)\big) & = T\big(\big[T^{|\varphi^n(2)|-2}(1), 1/\beta\big)\big)\,, \\[1ex]
T^{2g_n+1}\big(\big[f_{2\varphi^n(2)}^{-1}(1), f_{2\varphi^n(1)}^{-1}(1)\big]\big) & = T\big(\big[1/\beta, T^{|\varphi^n(1)|-1}(1)\big]\big)\,, \\[1ex]
T^{2g_n+1}\big(\big(f_{2\varphi^n(1)}^{-1}(1), T^{|\varphi^n(1)|-1}(1)\big]\big) & = T^{g_n+1}\big(\big[T^{|\varphi^n(2)|-2}(1), 1/\beta\big)\big)\,, \\[1ex]
T^{2g_n+1}\big(\big[T^{|\varphi^n(2)|-2}(1), 1/\beta\big)\big) & = T^{g_n}\big(\big[T^{|\varphi^n(21)|-1}(1), T^{|\varphi^n(1)|-1}(1)\big)\big)\,.
\end{align*}
Here, $f_{\varphi^n(21)}(1) \le 1$ implies that $T^{|\varphi^n(21)|-1}(1)  = f_2^{-1} f_{\varphi^n(21)}(1) \ge 1/\beta$.
Hence, the discontinuities of $T^{2g_n+1}$ in the interior of $F_{n,g_n-1}$ are $1/\beta$,  $f_{2\varphi^n(2)}^{-1}(1)$ and $f_{2\varphi^n(1)}^{-1}(1)$, with
\[
T^{2g_n+2}\big(\big[1/\beta, f_{2\varphi^n(2)}^{-1}(1)\big)\big) = F_{n+1,0} \setminus \{1\}\,, \ T^{2g_n+1}\big(\big[f_{2\varphi^n(2)}^{-1}(1), f_{2\varphi^n(1)}^{-1}(1)\big]\big) = F_{n,0}\,.
\]
This gives $T^m(I) \supseteq F_{n+1,0} \setminus \{1\}$ for some $m \ge 0$.

\medskip\noindent
\textbf{Case~3:} even~$n \ge 2$, $\eta_{n+1} \le \beta < \gamma_n$.
In this case, we have $\beta^{2g_n} = \beta^{g_{n+1}-1} > \beta^{g_n-1} + \beta^{-1} \ge \beta + \beta^{-1} > 2$, $F_{n,g_n-1} = \big[T^{|\varphi^n(1)|-2}(1), T^{|\varphi^n(2)|-1}(1)\big]$, and
\begin{align*}
T^{2g_n}\big(\big[1/\beta, f_{2\varphi^n(1)}^{-1}(1)\big)\big) & = T^{g_n}\big(\big[T^{|\varphi^n(1)|-2}(1), 1/\beta\big)\big)\,, \\[1ex]
T^{2g_n}\big(\big(f_{11\varphi^n(1)}^{-1}(1), 1/\beta\big)\big) & = T^{g_n-1}\big(\big(T^{|\varphi^n(1)|-2}(1), 1/\beta\big)\big)\,, \\[1ex]
T^{2g_n}\big(\big[f_{2\varphi^n(1)}^{-1}(1), T^{|\varphi^n(2)|-1}(1)\big]\big) & = T^{g_n}\big(\big[1/\beta, T^{|\varphi^n(21)|-2}(1)\big]\big)\,, \\[1ex]
T^{2g_n}\big(\big[T^{|\varphi^n(1)|-2}(1), f_{11\varphi^n(1)}^{-1}(1)\big]\big) & = T^{g_n-1}\big(\big[1/\beta, T^{|\varphi^n(2)|-1}(1)\big]\big)\,.
\end{align*}
Here, $f_{2\varphi^n(1)}\, T^{|\varphi^n(2)|-1}(1) = f_{\varphi^n(21)}(1) \ge 1$ implies $f_{2\varphi^n(1)}^{-1}(1) \le T^{|\varphi^n(2)|-1}(1)$. 
The discontinuities of $T^{2g_n}$ in the interior of $F_{n,g_n-1}$ are at most $1/\beta$, $f_{2\varphi^n(1)}^{-1}(1)$ and $f_{11\varphi^n(1)}^{-1}(1)$, with
\begin{align*}
T^{g_n+2}\big(\big[1/\beta, f_{2\varphi^n(1)}^{-1}(1)\big)\big) & = F_{n,0} \setminus \{1\}\,, \\ T^{g_n+3}\big(\big(f_{11\varphi^n(1)}^{-1}(1), 1/\beta\big)\big) & = F_{n,0} \setminus \{1, T^{|\varphi^n(1)|}(1)\}\,,
\end{align*}
similarly to Case~1. 
We obtain that $T^m(I) \supseteq F_{n,0} \setminus \{1, T^{|\varphi^n(1)|}(1)\}$ for some $m \ge 0$.
Since $T^{g_n+1}\big(F_{n,0} \setminus \{1, T^{|\varphi^n(1)|}(1)\}\big) \supseteq F_{n,0} \setminus \{1\}$, we also have $T^m(I) \supseteq F_{n,0} \setminus \{1\}$ for some $m \ge 0$.

\medskip\noindent
\textbf{Case~4:} even~$n \ge 2$, $\gamma_{n+1} \le \beta < \eta_{n+1}$. 
Since $\beta^{g_{n+1}+1} \ge \beta + 1 > 2$, we look at $T^{2g_n+2}$ on~$F_{n,g_n-1}$.
We have $T^{2g_n+1}\big(\big[1/\beta, f_{2\varphi^n(1)}^{-1}(1)\big)\big) = F_{n,g_n-1} \setminus \{T^{|\varphi^n(1)|-2}(1)\}$ with $f_{2\varphi^n(1)}^{-1}(1) > T^{|\varphi^n(2)|-1}(1)$, thus $f_{2\varphi^n(1)}^{-1}(1)$ is outside of $F_{n,g_n-1}$, but $T^{2g_n+2}$ still has a discontinuity in $\big(1/\beta, T^{|\varphi^n(2)|-1}(1)\big)$.
More precisely,
\begin{align*}
T^{2g_n+2}\big(\big[1/\beta, f_{2\varphi^n(2)}^{-1}(1)\big]\big) & = T\big(\big[1/\beta, T^{|\varphi^n(2)|-1}(1)\big]\big)\,, \\[1ex]
T^{2g_n+2}\big(\big(f_{11\varphi^n(1)}^{-1}(1), 1/\beta\big)\big) & = T^{g_n+1}\big(\big(T^{|\varphi^n(1)|-2}(1), 1/\beta\big)\big)\,, \\[1ex]
T^{2g_n+2}\big(\big(f_{2\varphi^n(2)}^{-1}(1), T^{|\varphi^n(2)|-1}(1)\big]\big) & = T\big(\big[T^{|\varphi^n(2)|-2}(1), 1/\beta\big)\big)\,, \\[1ex]
T^{2g_n+2}\big(\big[T^{|\varphi^n(1)|-2}(1), f_{11\varphi^n(1)}^{-1}(1)\big]\big) & = T\big(\big[T^{|\varphi^n(1)|-2}(1), T^{|\varphi^n(21)|-2}(1)\big]\big)\,,
\end{align*}
with $T^{|\varphi^n(21)|-2}(1) = f_{11}^{-1} f_{\varphi^n(21)}(1) < 1/\beta$ since $f_{\varphi^n(21)}(1) < 1$.
Therefore, the discontinuities of $T^{2g_n+2}$ in the interior of $F_{n,g_n-1}$ are $1/\beta$, $f_{2\varphi^n(2)}^{-1}(1)$ and $f_{11\varphi^n(1)}^{-1}(1)$, with
\begin{align*}
T^{2g_n+3}\big(\big[1/\beta, f_{2\varphi^n(2)}^{-1}(1)\big]\big) & = F_{n+1,0}\,, \\ T^{g_n+3}\big(\big(f_{11\varphi^n(1)}^{-1}(1), 1/\beta\big)\big) & = F_{n,0} \setminus \{1, T^{|\varphi^n(1)|}(1)\}\,.
\end{align*}
This gives $T^m(I) \supseteq F_{n+1,0} \setminus \{1\}$ for some $m \ge 0$.
\ep\medbreak

\begin{lemma} \label{l:covering}
Let $1 < \beta < \gamma_1$ and $n \ge 1$ such that $\eta_{n+1} < \beta < \eta_n$ with odd $n \ge 1$, or $\eta_{n+1} \le \beta \le \eta_n$ with even $n \ge 2$. 
Then we have
\[
T_{-\beta}^{g_n(2g_n+(-1)^n)}\big(F_{n,0}(\beta) \setminus \{1\}\big) \supseteq F(\beta) \setminus \big\{T_{-\beta}^k(1) \mid 0 \le k \le g_n(2g_n+(-1)^n)\big\}\,.
\]
\end{lemma}

\proc{Proof.}
Let $1 < \beta < \gamma_1$ and $n$ as in the statement of the lemma.
By Corollary~\ref{c:TFG}, we have $F(\beta) = \bigcup_{I \in \mathcal{F}_n(\beta)} I$ if $\beta < \gamma_n$, and $F(\beta) = \bigcup_{I \in \mathcal{F}_{n-1}(\beta)} I$ if $\beta \ge \gamma_n$.
In the latter case, $f_{\varphi^n(1)} \le f_{\varphi^{n-1}(1)}$ implies that $F_{n-1,0} = F_{n,|\varphi^{n-1}(1)|} \cup F_{n,0}$, and we obtain, similarly to the proof of Proposition~\ref{p:partition}, that $\bigcup_{I \in \mathcal{F}_{n-1}(\beta)} I = \bigcup_{I \in \mathcal{F}_n(\beta)} I$.
Therefore, we always have $F(\beta) = \bigcup_{I \in \mathcal{F}_n(\beta)} I$.

Let first $n \ge 1$ be odd.
For simplicity, we omit points in $\{T^k(1) \mid k \ge 0\}$ in the following statements.
We have $T^{g_n+1}(F_{n,0}) = F_{n+1,0} \cup F_{n,1}$ and
\[
T^{g_n-1}(F_{n+1,0}) = \big[T^{|\varphi^n(21)|-1}(1),\, T^{|\varphi^n(1)|-1}(1)\big]\,.
\]
Since $\beta > \eta_{n+1}$, we have $T^{|\varphi^n(21)|-1}(1) < 1/\beta$, thus $T^{g_n}(F_{n+1,0}) \supseteq F_{n,0}$, and
\[
T^{2g_n+1}(F_{n,0}) \supseteq F_{n,0} \cup F_{n,1}\,.
\]
Inductively, we get $T^{(g_n-1)(2g_n+1)}(F_{n,0}) \supseteq \bigcup_{k=0}^{g_n-1} F_{n,k}$.
Then,  $T(F_{n,g_n-1}) = F_{n,g_n} \cup F_{n,0}$ yields that
\[
T^{g_n(2g_n-1)}(F_{n,0}) = T^{(g_n-1)(2g_n+1)+1}(F_{n,0}) = \bigcup_{k=0}^{g_n} F_{n,k} = F(\beta)\,.
\]
Here, we have only omitted points in $\big\{T^k(1) \mid 0 \le k \le g_n(2g_n-1)\big\}$.

Let now $n \ge 2$ be even.
Omitting $\big\{T^k(1) \mid k \ge 0\big\}$, we have $T^{g_n}(F_{n,0}) = F_{n,g_n} \cup F_{n+1,0}$~and
\[
T^{g_n-1}(F_{n+1,0}) = \big[T^{|\varphi^n(1)|-2}(1),\, T^{|\varphi^n(21)|-2}(1)\big]\,.
\]
Since $\beta \ge \eta_{n+1}$, we have $T^{|\varphi^n(21)|-2}(1) \ge 1/\beta$, thus $T^{g_n+1}(F_{n+1,0}) \supseteq F_{n,0}$, and 
\[
T^{2g_n+2}(F_{n,0}) \supseteq F_{n,0} \cup F_{n,1}\,.
\]
As above, we get $T^{(g_n-1)(2g_n+2)}(F_{n,0}) \supseteq \bigcup_{k=0}^{g_n-1} F_{n,k}$.
Since $T^{g_n+2}(F_{n,g_n-1}) \supseteq F_{n,g_n} \cup F_{n,0}$ and $T^{g_n+2}(F_{n,k}) \supseteq F_{n,k+1}$ for $0 \le k < g_n$, we obtain that
\[
T^{g_n(2g_n+1)}(F_{n,0}) = T^{(g_n-1)(2g_n+2)+g_n+2}(F_{n,0}) = \bigcup_{k=0}^{g_n} F_{n,k} = F(\beta)\,,
\]
up to some points in $\big\{T^k(1) \mid 0 \le k \le g_n(2g_n+1)\big\}$.
\ep\medbreak

\section{Proofs of the main results} \label{sec:proofs-main-results}

Now we are ready to prove Theorems \ref{t:gaps} and
\ref{t:main-mixing}.

\proc{Proof of Theorem~\ref{t:main-mixing}.}
For $1 < \beta < \gamma_1$, Lemma~\ref{l:cover-F00} and~\ref{l:covering} show that $T_{-\beta}$ is locally eventually onto on $F(\beta) \setminus \{T_{-\beta}^k(1) \mid 0 \le k \le g_n(2g_n+(-1)^n)\}$ for some $n \ge 1$.
For $\beta \ge \gamma_1$, this was shown in \cite[Proposition~8]{Gora07}.
The iterated $T_{-\beta}$-images of any interval of positive length clearly contain~$1/\beta$, thus they also contain $T_{-\beta}^k(1)$ for any $k \ge 0$.
Therefore, $T_{-\beta}$ is locally eventually onto on $F(\beta) = (0,1] \setminus G(\beta)$.
The other statements of the theorem are in Corollary~\ref{c:TFG} and Proposition~\ref{p:gap2}. 
\ep\medbreak

\proc{Proof of Theorem~\ref{t:gaps}.}
By Proposition~\ref{p:gap}, all the intervals in $\mathcal{G}_n(\beta)$ are gaps. 
By Theorem~\ref{t:main-mixing}, they are the only gaps. 
Calculating directly the numbers of the intervals of $\mathcal{G}_n(\beta)$ by the definition, we complete the proof.
\ep\medbreak

Theorem~\ref{t:limit} was proved in Section~\ref{sec:polyn-expans}, thus it only remains to prove Theorem~\ref{t:perron}. 
Recall that a~number $\beta > 1$ is an Yrrap number if $V_{-\beta} = \{T_{-\beta}^n(1) \mid n \ge 0\}$ is a finite set.

\proc{Proof of Theorem~\ref{t:perron}.}
Let $V'_{-\beta}$ be the set of numbers in $V_{-\beta}$ which are not right endpoints of gaps, and $J_x = (\max\{y \in V_{-\beta} \cup \{0\} \mid y < x\}, x)$ for each $x \in V'_{-\beta}$.
Then $\{J_x \mid x \in V'_{-\beta}\}$ is a partition of $F(\beta)$ up to finitely many points, and $T_{-\beta}(J_x)$ is a union of intervals $J_y$, $y \in V'_{-\beta}$, for each $x \in V'_{-\beta}$.
Define a matrix
\[
M_{-\beta} = (m_{x,y})_{x,y \in V'_{-\beta}} \quad \mbox{with} \quad m_{x,y} = \frac{\beta\, \lambda(J_x \cap T_{-\beta}^{-1}(J_y))}{\lambda(J_y)}\,,
\]
i.e., $m_{x,y}$ is the number of times that $J_y$ is contained in $T_{-\beta}(J_x)$.
Since
\[
\sum_{y \in V'_{-\beta}} m_{x,y}\, \lambda(J_y) = \beta\, \lambda(J_x)
\]
for all $x \in V'_{-\beta}$, $(\lambda(J_x))_{x \in V'_{-\beta}}$ is a positive eigenvector of $M_{-\beta}$ to the eigenvalue~$\beta$.
Since $T_{-\beta}$ is locally eventually onto on~$F(\beta)$, the matrix $M_{-\beta}$ is primitive, thus $\beta$ is the Perron--Frobenius eigenvalue of~$M_{-\beta}$, which is a Perron number.
\ep\medbreak

We finally show that the set of Parry numbers and the set of Yrrap numbers do not include each other.
For the definition of a Parry number, it is convenient to extend the domain of $T_\beta$ to $[0,1]$ by setting $T_\beta(1) = \beta - \lfloor\beta\rfloor$. 
Then $\beta>1$ is a Parry number if and only if $(T_\beta^n(1))_{n\ge0}$ is eventually periodic. 
We know that all Pisot numbers are both Parry numbers \cite{Bertrand77,Schmidt80} and Yrrap numbers~\cite{Frougny-Lai09}.
Therefore, the symmetric difference between the set of Parry numbers and the set of Yrrap numbers can only contain Perron numbers that are not Pisot numbers.

\begin{prop}
Let $\beta > 1$ with $\beta^4  = \beta + 1$, i.e., $\beta \approx 1.2207$.
Then $T_{-\beta}^{10}(1) = T_{-\beta}^5(1)$, and $(T_\beta^n(1))_{n\ge0}$ is aperiodic.
\end{prop}

\proc{Proof.}
It is easy to check that $T_{-\beta}^{10}(1) = T_{-\beta}^5(1)$.

On the other hand, let $\alpha$ be an algebraic conjugate of $\beta$ satisfying $|\alpha| > 1$, i.e., $\alpha \approx -0.2481 \pm 1.034 i$, and let $\sigma:\, \mathbb{Q}(\beta) \to \mathbb{Q}(\alpha)$ be the field homomorphism defined by $\sigma(\beta) = \alpha$.
Since $T_{\beta}(x) \in \{\beta x, \beta x - 1\}$, we have $|\sigma(T_\beta(x))| \ge |\alpha|\, |\sigma(x)|  - 1$ for all $x \in [0,1] \cap \mathbb{Q}(\beta)$, thus $|\sigma(T_\beta(x))| > |\sigma(x)|$ if $|\sigma(x)| > 1 /(|\alpha|-1)$.
We have $T_\beta^{35}(1) = \beta^{35} - \beta^{34} - \beta^{26} - \beta^{13} - \beta^4$, and one can check that $|\sigma(T_\beta^{35}(1))| > 1/(|\alpha|-1)$.
This implies that $(|\sigma(T_\beta^n(1))|)_{n\ge35}$ is a strictly increasing sequence, hence $(T_\beta^n(1))_{n\ge0}$ is aperiodic.
\ep\medbreak

\begin{prop}
Let $\beta > 1$ with $\beta^7  = \beta^6 + 1$, i.e., $\beta \approx 1.2254$.
Then $T_\beta^7(1) = 0$, and $(T_{-\beta}^n(1))_{n\ge0}$ is aperiodic.
\end{prop}

\proc{Proof.}
We have $T_\beta^7(1) = \beta^7 - \beta^6 - 1 = 0$.

Let $\alpha$ be an algebraic conjugate of $\beta$ satisfying $|\alpha| > 1$, i.e., $\alpha \approx 0.7802 \pm 0.7053i$, and let $\sigma:\, \mathbb{Q}(\beta) \to \mathbb{Q}(\alpha)$ be the field homomorphism defined by $\sigma(\beta) = \alpha$.
We have $T_{-\beta}(x) \in \{-\beta x + 1, -\beta x + 2\}$, thus $|\sigma(T_{-\beta}(x))| \ge |\alpha|\, |\sigma(x)|  - 2$ for all $x \in (0,1] \cap \mathbb{Q}(\beta)$.
It is more convenient to consider $\phi(x)$ with $\phi$ as in the Introduction, since $\phi(T_{-\beta}(x)) \in \{-\beta \phi(x), -\beta \phi(x) - 1\}$ and $|\sigma(\phi(T_{-\beta}(x)))| \ge |\alpha|\, |\sigma(\phi(x))|  - 1$ for all $x \in (0,1] \cap \mathbb{Q}(\beta)$.
Now one can check that $|\sigma(\phi(T_{-\beta}^{53}(1)))| > 1/(|\alpha|-1)$, which implies that $(|\sigma(\phi(T_{-\beta}^n(1)))|)_{n\ge53}$ is a strictly increasing sequence, hence $(T_{-\beta}^n(1))_{n\ge0}$ is aperiodic.
\ep\medbreak

\section*{Acknowledgements}
Lingmin Liao was partially supported by NSFC10901124.

\bibliographystyle{amsplain}
\bibliography{minusbeta}
\end{document}